\documentclass[12pt,oneside,a4paper]{article}
\usepackage{fancyhdr, ifpdf}
\usepackage{syntonly}
\usepackage{geometry}
\usepackage{makeidx}
\makeindex
\usepackage[english]{babel}
\newtheorem{definition}{Definition}
\newtheorem{theorem}{Theorem}

\newtheorem{proposition}{Proposition}

\usepackage{makeidx}
\makeindex

\usepackage{amsmath}
\usepackage{amssymb}
\usepackage{amsmath}
\usepackage{amsfonts}
\usepackage{amssymb}
\usepackage{diagxy}

\title{Martin-L\" of randomness, invariant measures and  countable homogeneous structures}

\author{Willem L. Fouch\'{e}\\
\it Department of Decision Sciences,\\\it University of South
Africa, PO Box 392, 0003 Pretoria, South
Africa\\fouchwl@gmail.com}

\date{}

\begin{document}
 \maketitle

\pagenumbering{arabic}

\abstract\small{We  use  ideas from topological dynamics (amenability), combinatorics (structural Ramsey theory) and model theory (Fra\" {i}ss\' e limits)  to study  closed amenable subgroups $G$ of the symmetric group $S_\infty$ of a countable set, where $S_\infty$ has the topology of pointwise convergence. We construct $G$-invariant measures on  the universal minimal flows associated with these  groups $G$ in, moreover,  an algorithmic manner. This leads  to an identification of the generic elements, in the sense of being  Martin-L\" of random, of these flows with respect to the constructed  invariant measures. Along these lines we study the  random elements of $S_\infty$, which are permutations that transform recursively presented universal structures into such structures which are Martin-L\" of random.}

\paragraph{Keywords:} Martin-L\" of randomness, topological dynamics, amenable groups, \newline Fra\"iss\'e limits, Ramsey theory.

\newcommand{\qed}{\nobreak \ifvmode \relax \else
      \ifdim\lastskip<1.5em \hskip-\lastskip
      \hskip1.5em plus0em minus0.5em \fi \nobreak
      \vrule height0.30em width0.4em depth0.25em\fi}



                                                

\section{Introduction}

During the past four decades there has been a vigorous development on the interplay between combinatorics and algorithmics on the one hand, and the dynamical properties of topological groups, on the other. (See, for example,  Sarnak \cite{Sarnak},  as well as Lubotzky, Phillips and Sarnak  \cite{LubotzkyPhillipsSarnak1, LubotzkyPhillipsSarnak2}.)

Notions from the theory of topological transformation groups such as amenability or Kazhdan's property (T), play a central r\^ole in the construction of expander graphs (as expounded for example in \cite{Sarnak}), which in turn are central to dealing with deterministic error amplification for the complexity class $\mathbf{RP}$ (randomized polynomial time algorithms). The aim of this very active area of research still is to minimise  the number of random bits (as generated, for example,  by some physical artifact) which might be required for the execution of a probabilistic computation. These results also have highly nontrivial implications for the design of quantum circuits. (See, for example, Harrow, Recht and Chuang  \cite{HarrowRechtChuang}.)  We are dealing here with  essential instances of the deep problem of {\it derandomization}. 

The focus of this paper is on the dual problem of understanding the symmetries that transform a recursively presented  universal structure, which in this paper is a Fra\" {i}ss\' e limit of finite first order structures, to a copy of such a structure which is Martin-L\" of random relative to an $S_\infty$-invariant measure on the class of all universal structures of the given type. Here $S_\infty$ is  the symmetric group of a countable set, with the pointwise convergence topology.

It was shown by Kechris, Pestov and Todorcevic \cite{Kechrisetal}, that deep  results in combinatorics (structural Ramsey theory) can be interpreted as statements on the dynamical properties of closed subgroups of $S_\infty$.  They identified the so-called extremely amenable subgroups $G$ of $S_\infty$  in terms of the Ramsey properties of some classes $\mathcal{F}$ of finite first-order structures, such that $G$ is the symmetry group of the Fra\" {i}ss\' e limit $\mathbb{F}$ of the class $\mathcal{F}$. They also showed  how one can utilise the Ramsey properties of some classes $\mathcal{F}$ to identify the universal minimal flow $\mathcal{U}$  of  the automorphism group  $G$ of the Fra\" {i}ss\' e limit $\mathbb{F}$ of $\mathcal{F}$.  It is a remarkable fact  that for many interesting classes  $\mathcal{F}$ (total orders, graphs, posets, ranked diagrams, \dots) the space $\mathcal{U}$ can be embedded into a Baire space of the form $\{0,1\}^{\mathbb{N}^k}$ which renders $\mathcal{U}$ accessible to an effective study of randomness relative to computable measures on $\{0,1\}^{\mathbb{N}^k}$  which are supported by $\mathcal{U}$. 

In this paper we initiate a study of computable group invariant measures $\nu$ on the spaces $\mathcal{U}$ and identify the elements of $\mathcal{U}$ which are Martin-L\" of random with respect to the measures $\nu$. We shall study the elements of $S_\infty$ (called randomizers in this paper) that transform universal recursive objects in $\mathcal{U}$ to such objects which are Martin-L\" of random with respect to $\nu$. We shall also look at the recent paper of Petrov and Vershik \cite{PetrovVershik1} from the viewpoint of randomizing recursive universal structures.

\section{Preliminaries on amenable groups}
\label{sec:amenable groups}
   Issues involving the uniqueness of Lebesgue measure, including the Banach-Tarski paradox, led to questions of how  a group $G$ acting on a set $X$ determines the structure (and
number) of $G$--invariant finitely additive probability measures on $X$ ($G$-invariant means). For example, Banach (1923) \cite{Banach} showed that there was more
than one rotation invariant finitely additive probability measure on the Lebesgue measurable subsets of $S^1$. Subsequently, Sullivan (1981) \cite{Sullivan} and independently
Margulis (1980, 1981) \cite{Margulis1} and \cite{Margulis2} (for $n \geq 4$) and Drinfield (1984) \cite{Drinfield} (for $n = 2,3$) showed
that Lebesgue measure is the unique finitely additive rotation invariant
measure on the Lebesgue measurable subsets of $S^n$. (It is crucial to consider the Lebesgue measurable sets here. The corresponding problem for Borel--measureable sets is still open.) It is probably fair to say that the entire development of amenable groups arose from this interest.

Let $G$ be a topological group and $X$ a compact Hausdorff space. A dynamical
system $(X,G)$ (or a $G$-flow on $X$) is given by a jointly continuous action of $G$ on $X$. If $(Y, G)$ is
a second dynamical system, then a $G$-morphism  $\pi : (X, G) \rightarrow  (Y, G) $
is a continuous mapping $\pi:X \rightarrow Y$ which intertwines the  $G$-actions, i.e.,the diagram

\[\bfig
\square[G \times X`X`G\times Y`Y; \alpha`1\times\pi`\pi`\beta  ]
   \efig ,\]    
commutes with $\alpha,\beta$ being the group actions.
 
 An isomorphism is a bijective homomorphism. A subflow of $(X,G)$ is a $G$-flow on a compact subspace $Y$ of $X$ with the action the restricted to the action of $G$ on $X$ to the action on $Y$. The dynamical
system $(X, G)$ is said to be \emph{minimal} if every $G$-orbit is dense in $X$. Equivalently, a $G$-flow is minimal if it has no proper subflows. Every dynamical system has a minimal subflow (Zorn).                         

The following fact, first proven by Ellis (1949) \cite{Ellis}, is central to the theory of dynamical systems:
\begin{theorem}
Let $G$ be a Hausdorff topological group. There exists, up to $G$-iso-\newline morphism, a unique minimal dynamical system, denoted by $(M (G),G)$, such that   for every minimal dynamical system $(X, G)$ there exists a $G$--epimorphism $\pi : (M, G) \rightarrow
(X, G)$, and any two such universal systems are isomorphic.
\end{theorem}
The flow $(M (G),G)$ is called the {\it universal minimal flow} of $G$. We next introduce the notion of amenable groups.
\begin{definition}
A topological group $G$ is {\it amenable} if, whenever $X$ is a non-empty compact Hausdorff space and $\pi$ is a continuous action of $G$ on $X$, then there is a $G$--invariant Borel probability measure on $X$.
\end{definition}
This means that, for every $G$-flow on a compact space $X$,  there is a measure $\nu$ on the Borel algebra of $X$, such that, $ \nu(X)=1$ and, for every $g \in G$ and Borel subset $U$ of $X$,
$$\nu(gU)=\nu(U).$$
It follows that $G$ is amenable iff its universal minimal flow $M(G)$ has a $G$-invariant probability measure. Indeed, let $\nu$ be an invariant measure on $M(G)$. Consider any $G$-flow on some compact Haussdorf space $X$. By Zorn's lemma there is a minimal subflow $Y$ and a $G$-embedding $i$ of $Y$ into $X$. Therefore, there  are $G$-morphisms
\[\bfig
    \morphism (0,0) <500,0>[M(G)`Y;\pi]
     \morphism (500,0)/{>->}/<450,0>[Y`X; i]
\efig.\]
Let $\rho$ be the pushout measure of $\nu$ under $i\pi$.  In other words, for every Borel subset $A$ of $X$, we set
\[\rho(A) = \nu(\pi^{-1}i^{-1}A).\]
Then $\rho$ is an invariant measure on $X$. The converse is trivial, since $M(G)$ is a compact $G$-flow. 

If $G$ is compact, then it is well-known that $M(G) \simeq G$, topologically, and the Haar measure is therefore the (unique) invariant measure on $M(G)$. In particular, a compact Hausdorff group is amenable. It can be shown that a discrete group $G$ is amenable (see, for example \cite{Fremlin}) iff it admits a $G$--invariant finitely additive probability measure on all the subsets of $G$. It is a classical fact that the free group $\mathbb{F}_2$ on two elements is not amenable; this result played an important r\^ ole in understanding the Banach-Tarski paradox.

We  make frequent use of the following:
\begin{theorem}
 Let $G$ be a topological group and suppose there is a dense subset of $G$ such that every finite subset of $G$ is included in an amenable subgroup of $G$, then $G$ too is amenable.
\label{th:fremlin}
\end{theorem}

A proof of this theorem can be found in Section $449C$ of \cite{Fremlin}. 

For notational convenience, we shall assume in the sequel that the elements of $S_\infty$  are permutations of the natural numbers $\mathbb{N}$.   There is natural embedding of $S_\infty$ into the space $\mathbb{N}^\mathbb{N}$. We give $\mathbb{N}^\mathbb{N}$ the product topology (with $\mathbb{N}$ having the discrete topology). We topologise $S_\infty$with the inherited topology from $\mathbb{N}^\mathbb{N}$. For obvious  reasons, this is called the topology of pointwise convergence. 

For a natural number $n$ and two $n$-tuples $(a_1, \ldots, a_n)$ and $(b_1, \ldots, b_n)$ of natural numbers, each $n$-tuple having $n$ distinct elements, we denote by
\begin{equation*}
\begin{bmatrix}

a_1  &a_2 &  \cdots &a_ n\\
b_1 &b_2 & \cdots & b_n 
\end{bmatrix}
\end{equation*}
the set of elements of $S_\infty$ that map $a_i$ to $b_i$ for all $i=1,\ldots,n$. Each of these sets is open in $S_\infty$ and all these (cylindrical) sets constitute  a basis for the topology on $S_\infty$.

Since any finite group (with the discrete topology), being compact,  is amenable, it follows from Theorem \ref{th:fremlin}  that if a group contains a dense locally finite subgroup, then it is amenable.
(Recall that a  group is locally finite if every finite subset of the group generates a finite group.)
 In particular, any locally finite group is amenable. Let $S_0$ be the subgroup of $S_\infty$ consisting of all the permutations that fix all but a finite set of elements. Then $S_0$ is dense in $S_\infty$ and  the former is obviously locally finite. It follows that $S_\infty$ is an amenable group.

In contrast to locally compact groups, it is readily seen that not all closed subgroups of $S_\infty$ are amenable groups. Indeed, a  Caley embedding would embed any countable group as a closed discrete subgroup of $S_\infty$. Indeed, let $G$ be a countable group. The Caley embedding $c$  is defined to be
\[ c: G \longrightarrow S_\infty,\]
given by
\[ \sigma \mapsto (\sigma),\]
where $(\sigma)$ is the permutation on $\mathbb{N}$ given by $x \mapsto \sigma x$, for all $x \in \mathbb{N}$. It is clear that the image of $G$ under $c$ inherits from $S_\infty$ the discrete topology. Moreover, $c(G)$ is easily seen to be closed in $S_\infty$. (In fact, any discrete subgroup in a Hausdorff topological group is necessarily closed!)

In particular, the non-amenable group $\mathbb{F}_2$ can be embedded as a discrete and closed subgroup of $S_\infty$.
 
We next discuss a beautiful example of Pestov \cite{Pestov1}: Let $\eta$ be the countable linear order, called the rational order, which is order isomorphic to the standard ordering of the rational numbers. Let $A$ be the automorphism group of the ordering $\eta$. Then $A$ is a closed subgroup of $S_\infty$. As has already been mentioned,  the free group on two elements $\mathbb{F}_2$ with the discrete topology is not an amenable group.  Let $(D, <)$ be a countable linearly ordered skew field. Then the multiplicative group $D^\times$ of $D$ acts via left-multiplication on the total order $<$. 

 However, as shown by Neumann (1949) \cite{Neumann}, the order $<$ is isomorphic to $\eta$. In this way one embed $D^\times$ as a closed and discrete subgroup of $A$. 

Finally, Neumann also showed that $\mathbb{F}_2$ can be embedded into $D^\times$. In this way we get a closed embedding of $\mathbb{F}_2$ in $A$. 

This construction is extremely interesting, for the group $A$ is in fact extremely amenable (no pun intended).  This means that every continuous action of $A$ on a compact metric space has a fixed point.

 In \cite{ BhattacharjeeMacpherson}  Bhattacharjee and Macpherson have shown that there is a locally finite dense  subgroup of the symmetry group of the random graph, thus allowing us to infer that the symmetry group of the well-known random graph is amenable. We shall later show  how one can algorithmically construct an invariant mean on the universal minimal flow of the symmetry group of the random graph. 

It was recently shown by Kechris and Soli\v c \cite{KechrisSokic} that the symmetry group of the Fra\" {i}ss\' e  limit (see the following section for definitions) of finite posets   is not amenable. The Ramsey properties of finite posets (see, for example, \cite{Fouchposet}) play an important role in their argument. These results will later be used in this paper to show that the set of linear extensions of the Fra\" {i}ss\' e  limit  of finite posets  are all, in a definite sense, nonrandom, at least from the point of view of algorithmic randomness.

\section{Fra\" {i}ss\' e limits and their recursive representations}
\label{sec:fraisse limits}
In the sequel, ${\cal L}$ will stand for the 
signature of 
a relational 
structure. Moreover, ${\cal L}$ will always be finite and 
the arities 
of the relational 
symbols will all be $\geq 1$.
 The definitions that follow were 
introduced by 
Fra\"\i ss\'e 
in 1954.

The {\em age} of an ${\cal L}$-structure $X$, written $\mbox{\tt Age}( X )$, is 
the class of all finite 
${\cal L}$-structures (defined on finite ordinals) 
which can 
be embedded as ${\cal 
L}$-structures into $X$. The structure $X$ is {\em 
homogeneous} (
some authors say {\em
ultrahomogeneous} ) if, given any isomorphism $f: A \rightarrow
B$ between 
finite substructures of 
$X$, there is an automorphism $g$ of $X$ whose restriction 
to $A$ is 
$f$. 
A class {\bf K} of finite ${\cal 
L}$-structures has the 
{\em amalgamation 
property} if, for structures $A, B_1, B_2$ in {\bf K} and 
embeddings 
$f_i : A \rightarrow B_i$ 
($i=1,2$) there is a structure $C$ in {\bf K} and there are 
embeddings 
$g_i : B_i \rightarrow C$ 
($i=1,2$), such that the following diagram  commutes: 

   \[\bfig
   \morphism(0,0)|l|<250,-250>[A` B_2;f_2] 
   \morphism(0,0)<250,250>[A`B_1;f_1] 
   \morphism(250,250)/{-->}/<250,-250>[B_1`C;g_1]
    \morphism(250,-250)|r|/{-->}/<250,250>[B_2`C ;g_2] 
      \efig.\]
Suppose {\bf K} is a countable class of finite ${\cal 
L}$-structures, 
the domains of which are 
finite ordinals such that 
\begin{enumerate}
\item if $A$ is a finite ${\cal L}$-structure defined on 
some finite 
ordinal,  if $B \in 
\mbox{\bf K}$ and if there is an embedding of $A$ into $B$, 
then $A \in 
\mbox{\bf K}$;
\item the class {\bf K} has the amalgamation property.
\end{enumerate}
Then, Fra\"\i ss\'e showed that there is a countable 
homogeneous 
structure $X$ such that 
$\mbox{\tt Age}(X) = \mbox{\bf K}$. Moreover, $X$ is unique up to 
isomorphism. 
The (essentially) unique $X$ is called 
the {\em Fra\"\i ss\'e limit}  $\mathbb{K}$ of {\bf K}. Note
that, 
conversely, the age {\bf K} of 
a countable homogeneous structure has properties (1) and (2). We shall frequently call a countable structure which is isomorphic to a Fra\" {i}ss\' e limit a {\em universal structure}.

 A {\em recursive representation} of a countably infinite 
${\cal 
L}$-structure 
$X$ is a bijection $\phi: $X$ \rightarrow \mathbb{N}$ such that, for 
each $R \in 
{\cal L}$,
if the arity of $R$ is $n$, then the relation $R^{\phi}$ 
defined on 
$\mathbb{N}^n$ by
$$R^{\phi} \left( x_1, x_2, \ldots, x_n \right) 
\Longleftrightarrow R 
\left( \phi^{-1}(x_1), 
\ldots, \phi^{-1}(x_n) \right),$$
is recursive. If we identify the underlying set of $X$ with 
$\mathbb{N}$ 
via $\phi$ and each $R$ 
with $R^{\phi}$, we call the resulting structure a {\em 
recursive ${\cal 
L}$-structure} on $\mathbb{N}$ and we say it has a recursive representation on $\mathbb{N}$.

If $X$ is countable and homogeneous and if $\mbox{\tt Age} (X)$ has an 
enumeration $A_0$, $A_1$, $A_2$, 
$\ldots $,
possibly with repetition, with the property that there is a 
recursive 
procedure that 
yields, for each $i \in \mathbb{N}$, and $R \in {\cal L}$, the 
underlying set 
$A(i)$ of 
$A_i$ together with the interpretation of $R$ in $A(i)$, 
then we call 
$\left( A_i : i  \in \mathbb{N}  \right)$ a recursive enumeration of 
$\mbox{\tt Age} (X)$ . It 
follows from the construction of Fra\"\i ss\'e limits from 
their ages, 
that one can construct a recursive 
representation of $X$ 
from a recursive enumeration of its age. (Conversely, 
it is 
trivial to 
derive a recursive enumeration of $\mbox{\tt Age} (X)$ from a recursive 
representation of $X$.) 
It is therefore not difficult to find recursive 
representations for 
Fra\"\i ss\'e 
limits of classes {\bf K} from recursive enumerations of 
their ages.
\begin{theorem}
Suppose $\mathbb{C}$ and $\mathbb{D}$ are countable recursively represented   $\cal{L}$-structures on $\mathbb{N}$ with the same age. Suppose that they are both homogeneous. Then there is a recursive isomorphism from $\mathbb{C}$ to $\mathbb{D}$.
\label{th:efffraisse}
\end{theorem}
Proof. The model-theoretic back-and-forth argument as  discussed, for example, on pp 161-162 of Hodges \cite{Hodges} is constructive relative to the recursive representations of the homogeneous structures $\mathbb{C}$ and $\mathbb{D}$.

Let $\mathcal{K}$ be a Fra\"\i ss\'e class of finite structures. We say that $\mathcal{K}$ has the
{\it Hrushovski property} if for any $A$ in $\mathcal{K}$ there is $B$ in $\mathcal{K}$ containing $B$ such
that any partial automorphism of $A$ extends to an automorphism of $B$. The terminology, due to Kechris and Rosendal \cite{KechrisRosendal}, is probably inspired by the result by Hrushovski \cite{Hrushovski} who established that the class of finite graphs has this property. In \cite{KechrisRosendal} it is shown that a Fra\"\i ss\'e class of finite structures $\mathcal{K}$
has the Hrushovski property iff the automorphism group $G$ of the Fra\"\i ss\'e limit of $\mathcal{K}$ is compactly approximable, i.e., there is a increasing sequence $K_n$ of compact subgroups whose union is dense in the automorphism group. Since a compact group is amenable, it follows from Theorem \ref{th:fremlin} that the automorphism group of the Fra\"\i ss\'e limit
of a Fra\"\i ss\'e class with the  Hrushovski property will be amenable.

\section{Martin-L\"of random countable orders}
\label{sec:random countable orders}


Let $S_\infty$ be the group of permutations  of a countable set, which, without loss of generality, we may take to be  $\mathbb{N}$. We place on $S_\infty$ the pointwise convergence topolopy.

Let $(\mathbb{N}\times\mathbb{N})_{\neq}$ denote the set of ordered pairs $(i,j)$ of natural numbers with $i \neq j$. Write $\mathcal{M}$ for the set of total orders on  $\mathbb{N}$. We identify $\mathcal{M}$ with a subset of $\{0,1\}^{(\mathbb{N}\times\mathbb{N})_{\neq}}$ by identifying a total order  $<$  on $\mathbb{N}$ with the function $\xi:(\mathbb{N}\times\mathbb{N})_{\neq} \rightarrow \{0,1\}$ given by
$$\xi(x,y)=1 \Leftrightarrow x < y,\;\;x,y \in \mathbb{N}.$$ The total order associated with $\xi$ will be denoted by $<_\xi$.
We topologise $\mathcal{M}$ via the natural injection $\mathcal{M} \longrightarrow \{0,1\}^{(\mathbb{N}\times\mathbb{N})_{\neq}}$, where the (Baire) space $\{0,1\}^{(\mathbb{N}\times\mathbb{N})_\neq}$ has the product topology. As such $\mathcal{M}$ is a closed hence compact subspace of $\{0,1\}^{(\mathbb{N}\times\mathbb{N})_\neq}$.

 For each finite order $l$ on a subset of $\mathbb{N}$, write
$A_l$ for the space of total orders $\xi$ on $\mathbb{N}$ which extends $l$. Then the class $(A_l)$ with $l$ ranging over all the finite total orders on subsets of $\mathbb{N}$ is an open basis of neighbourhoods of  $\mathcal{M}$. The group $S_\infty$ acts continuously on $\mathcal{M}$ if, for $\xi \in \mathcal{M}$ and $\sigma \in S_\infty$, we define the total order $\sigma\xi$ by:
\[x <_{\sigma\xi} y \Longleftrightarrow \sigma^{-1}x <_\xi \sigma^{-1}y,\;x,y \in \mathbb{N}.\]

It follows from the classical (finitary) Ramsey theorem and the fundamental paper by Kechris, Pestov and Todorcevic (2005) that the  group action 
\[ S_\infty \times \mathcal{M} \longrightarrow \mathcal{M},\]
\[(\sigma, \xi) \mapsto \sigma\xi\]
is indeed isomorphic to the universal minimal $S_\infty$-flow.  The original  proof of this result can be found in the paper \cite{GlasnerWeiss} of Glasner and Weiss (2002). Their proof is also based on the classical Ramsey theorem.

Since $S_\infty$ is an amenable group, there is an $S_\infty$-invariant   (Borel regular) probability measure on $\mathcal{M}$. In fact, Glasner and Weiss (2002) \cite{GlasnerWeiss} showed  that there is {\it exactly one} such measure (i.e., the flow on $\mathcal{M}$ is uniquely ergodic). 
Their proof is based on an ergodic argument.

Let us denote this measure by $\mu$. The author believes this beautiful measure deserves to be called the {\it Glasner-Weiss measure}.

In this paper, we wish to understand the $\mu$-Martin-L\"of random elements of $\mathcal{M}$.  This is a viable project, for it will now be shown, with the benefit of hindsight, that the Glasner-Weiss measure can be computed and effectively constructed. (For more on the subject of Martin-L\"of randomness relative to computable measures, the reader is referred to the recent survey paper \cite{ML_randomness} together with the many references to be found there.)

 We write $\mathcal{M}_f$ for the set of finite total orders on some subset of $\mathbb{N}$. For $\ell \in \mathcal{M}_f$,  denote by $Z_\ell$ the set of $\xi \in \mathcal{M}$, such that $\xi$ is an extension of $\ell$. These  sets are the cylinder subsets of $\mathcal{M}$. Write $\mathcal{Z}_0$ for the class of events of the form $Z_\ell$ for some $\ell \in \mathcal{M}_f$ and $\mathcal{Z}$ for the algebra generated by $\mathcal{Z}_0$. Note that the $\sigma$-algebra generated by
 $\mathcal{Z}$ is exactly the Borel algebra on $\mathcal{M}$.

Clearly, for $\sigma \in S_\infty$ and $\ell \in \mathcal{M}_f$, 
\[ \sigma Z_\ell  = Z_{\sigma\ell}.\]
Furthermore, writing $S_\ell$ for the subgroup of $S_\infty$ each element of which  permutes the underlying set of $\ell$ and leaves all the other elements of $\mathbb{N}$ fixed, we have the following partition:
\[\mathcal{M} =\Cup_{\sigma \in S_\ell} \sigma Z_\ell.\]

Since $S_\infty$ is amenable  and $\mathcal{M}$ is compact, there is an $S_\infty$-invariant  probability measure $\nu$ on the Borel subsets of $\mathcal{M}$.  For such a measure, it must follow
\[ 1= \nu(\mathcal{M}) = \sum_{\sigma \in S_\ell} \nu(\sigma Z_\ell).\]

Since $\nu$ is $S_\infty$-invariant, it follows that 
\begin{equation}
\nu(Z_\ell) =\frac{1}{l!},
\label{eq:basic}
\end{equation}
where $l$ is the cardinality of the underlying set of $\ell$. Since a Borel measure on $\mathcal{M}$ is uniquely determined by its values on $\mathcal{Z}_0$, the uniqueness of the invariant mean $\nu$ has been established. In particular, $\nu$ is the Glasner-Weiss measure $\mu$.

For $Z \in \mathcal{Z}_0$ we write $Z^0$ for the complement of $Z$ and $Z^1$ for $Z$. (This peculiar notation has been designed in order to justify an algorithm for computing $\mu$ on $\mathcal{Z}$, as will become clear in the sequel.) Let  $(T_i)_{i \in \mathbb{N}}$ be any enumeration of  the algebra $\mathcal{Z}$ generated by $(Z_\ell)_{\ell \in \mathcal{M}_f}$ in such a way that one can effectively retrieve from  a given $i \in \mathbb{N}$,  the corresponding  $T_i$ as a finite union of sets $T$ of the form 
\begin{equation}
T = Z_{\ell_1}^ {\delta_1} \cap \ldots Z_{\ell_k} ^{\delta_k},
\label{eq:boolean}
\end {equation}
where each $\ell_i$ is in $\mathcal{M}_f$ and $\delta_i \in \{0,1\}$ for $i=1, \ldots, k$. We call any such enumeration a {\it recursive representation} of $\mathcal{Z}$.

The measure $\mu$ is computable in the following sense: 

\begin{theorem} 
Denote by $\mu$ the Glasner-Weiss measure on the Borel-algebra of $\mathcal{M}$. Let $\left(T_{i} : i <\omega\right)$ be a recursive representation of the algebra $\mathcal{Z}$ . 
 There is an effective procedure that yields, for $i,k \in \mathbb{N}$,  a binary rational $\beta_{k}$ such that
$$|\mu\left(T_{i}\right)-\beta_{k}|<2^{-k}.$$
\end{theorem}
Proof.  By the principle of inclusion and exclusion, it suffices to find a recursive procedure for computing the $\mu$-measure of all sets of the form (\ref{eq:boolean}).  For such an $T$, we call $n-\left(\delta_{1}+\ldots +\delta_{n}\right)$ the \emph{weight} of the representation.  If follows from (\ref{eq:basic})  that we can compute $\mu(T)$ up to arbitrary accuracy for all $T$ having a representation of weight 0. Now if $T$ has a representation of weight $f+1$, say, we can effectively write $T$ as $T'\cap T_{i}^{0}$ where $T_{i}\in \mathcal{Z}_{0}$ and $T'$ has a representation of weight $f$.  However,
$$\mu\left(T'\right)=\mu\left(T'\cap T_{i}^{0}\right)+\mu\left(T'\cap T_{i}\right),$$
so we can recursively call the procedure for computing any $\mu(G)$  with $G$ having a representation of weight $f$, to compute $\mu\left(T'\right)$ and $\mu\left(T'\cap T_{i}\right)$ up to accuracy $2^{-(k+1)}$ and hence to compute $\mu(T)$ up to accuracy $2^{-k}$.
\begin{definition}
 A set $A\subset \mathcal{M}$ is of \emph{constructive measure} $0$, if, for some recursive representation of  $\left(T_{i}:i \in \mathbb{N}\right)$ of $\mathcal{Z}$, there is a total recursive $\phi : \mathbb{N}^{2}\rightarrow \mathbb{N}$ such that
$$A\subset \bigcap_{n}\bigcup_{m}T_{\phi(n,m)}$$
and $\mu\left(\bigcup_{m}T_{\phi(n,m)}\right)$ converges effectively to $0$ as $n\rightarrow \infty$.
\end{definition}
\begin{definition}
A total order $\xi$ is said to be $\mu$-{\emph Martin-L\" of random}  if $\xi$ is in the complement of every subset $B$ of $\mathcal{M}$ of constructive measure $0$. 
\end{definition}

Write $ML_\mu \subset \mathcal{M}$ for the set  of  $\mu$-Martin-L\"of random total orders. Note that $\mu(ML_\mu)=1$. We now prove the following
\begin{theorem}
Write $\mathcal{Q}$ for the set of total orders on $\mathbb{N}$ which are isomorphic to the rational order $\eta$. Then
\[ML_\mu \subset \mathcal{Q}.\]
In particular,
\[\mu(\mathcal{Q})=1.\]
\label{th:mlcantor}
\end{theorem}

We first introduce a number of (standard) recursion-theoretic concepts and terminology: A sequence $\left(A_{n}\right)$ of sets in $\mathcal{Z}$ is said to be \emph {semirecursive} if for each $n$, the set $A_n$  is of the form $T_{\phi(n)}$ for some total recursive function $\phi : \omega \rightarrow \omega$ and some effective enumeration $\left(T_{i}\right)$ of $\mathcal{Z}$.  (Note that the sequence $\left(A_{n}^{c}\right)$, where $A_{n}^{c}$ is the complement of $A_{n}$, is also an $\mathcal{Z}-$semirecursive sequence.) In this case, we call the union $\bigcup_{n}A_{n} $ a $\sum_{1}^{0}$ set.  A set is a $\prod_{1}^{0}$ set if it is the complement of a $\sum_{1}^{0}$ set.  It is of the form $\bigcap_{n}A_{n}$, for some $\mathcal{Z}$-semirecursive sequence $\left(A_{n}\right)$. The proof of Theorem \ref{th:mlcantor} is based on the following
\begin{proposition}
If $A$ is a $\sum_{1}^{0}$ subset of $\mathcal{M}$  and if $\mu(A)=1$, then $ML_\mu$ is contained in $A$. In particular, if $B$ is a $\Pi_1^0$ subset of $\mathcal{M}$ that contains some element of $ML_\mu$, then $\mu(B) >0.$ 
\label{prop:kurtz}
\end{proposition}
Proof.  Let $P$ be the complement of $A$.  We can write $P=\bigcap_{n} P_{n}$ for some semirecursive sequence $\left(P_{n}\right)$.  By replacing each $P_{n}$ by $\bigcap_{k<n}P_{k}$, the sequence remains semirecursive but now $P_{n}$ is decreasing in $n$.  Since $\lim_{n\rightarrow \infty} \mu\left(P_{n}\right)=0$ it follows, from the computability of $\mu$,  that, for each $n$ one can effectively find some $f(n)$ such that $\mu\left(P_{f(n)}\right)<\left(n+1\right)^{-1}$.  Since $\mu\left(P_{i}\right)$ is decreasing in $i$, we then have that $\mu\left(P_{m}\right)<\left(n+1\right)^{-1}$ for all $m\geq f(n)$.  We conclude that $P$ is of constructive measure 0. \\
Proof of Theorem \ref{th:mlcantor}. In view of Cantor's well-known first-order characteristion of the rational order, it will suffice to proof the following: If $\xi \in ML_\mu$, then, for all natural numbers $n,m$ with $n \neq m$,

\begin{equation}
\exists _j \; n<_\xi j <_\xi  m; 
\label{eq:cantor1}
\end{equation}
moreover, for every natural number $n$
\begin{equation}
\exists_{j,k}\; j <_\xi n  <_\xi k.
\label{eq:cantor2}
\end{equation}
 Note that, for fixed $n,m$,  both  of these  predicates predicates are $\Sigma_1^0$ in $\xi$. It follows from Proposition \ref{prop:kurtz} that Theorem \ref{th:mlcantor} will follow once we have shown that the predicates (\ref{eq:cantor1}) and (\ref{eq:cantor2}) define events both of which have $\mu$-measures $1$.

Write $C_{n,m}$ for the complement of the event defined by the predicate (\ref{eq:cantor1}). Note that 
\[C_{n,m} = \bigcap_N\; C_{n,m}^N,\]
where
\[\xi \in  C_{n,m}^N \Leftrightarrow  \forall_{j \leq N} \; j \leq_\xi   n \vee  j\geq_\xi  m.\]

It follows from (\ref{eq:basic}) that, for $N > n,m$, we have

\[ \mu( C_{n,m}^N) =\frac{L}{N!},\]
where $L$ is the number of total orders $\ell$ on $\{1, \dots, N\}$ such that $n,m$ are adjacent with respect to the total order $\ell$. Clearly
\[L = 2(N-1) (N-2)!.\]
Consequently, 
\[ \mu(C_{n,m}) = \lim_{N \rightarrow \infty} \mu( C_{n,m}^N) =0.\]
Hence, if $\xi \in ML_\mu$, then $\xi$ satisfies (\ref{eq:cantor1}) for all $n,m$ with $n \neq m$. The proof that $\xi$ will satisfy (\ref{eq:cantor2}) for all $n$ is similar. This concludes the proof of Theorem \ref{th:mlcantor}.\\
{\bf Remark}: 
By using standard recursion-theoretic techniques it is not difficult to show that $ML_\mu$ is in fact $\Sigma_2^0$ with a complement being a set of constructive measure $0$. In other words, there is a {\it universal}  $\mu$-Martin-L\" of random test.

Let $\mathbb{P}=(P, \prec)$ be the universal countable poset.  This is, by definition, the Fra\"{i}ss\'e limit of the class of finite posets.  If $A$ and $B$ are subsets of $P$, then we write $A \prec B$, if for all $a \in A$ and $b \in B$, it is the case that $ a \prec b$.
The structure  $\mathbb{P}$ has the following first order characterisation:

If $A,B,Z$ are finite and pairwise disjoint subsets of the underlying set $P$ of $\mathbb{P}$ such that  $A \prec B$ and, for all $a \in A, b \in B$ and $z \in Z$:
\[ \neg z \prec a, \neg b \prec z,\]
then there is some element $x$ in $\mathbb{P}$ such that $A \prec x \prec B$ and $x$ is incomparable with all the elements of $Z$.

Let $\mathcal{P}^*$ be the class of all  structures of the form $(\pi, <)$, where $\pi$ is a finite poset (with underlying set some ordinal), and where $<$ is a total order on the underlying set of $\pi$ which is, moreover, a linear extension of $\pi$. Then  $\mathcal{P}^*$ is a {{F}ra\"{i}ss{\'e} class and its limit is of the form $\mathbb{P}^*=(\mathbb{P}, <)$, where $\mathbb{P}$ is the universal poset and $<$ is some linear extension of $\mathbb{P}$. We call the linear extension $<$ of $\mathbb{P}^*$ the {\it canonical} linear extension of $\mathbb{P}$. It readily follows from the fact that the structure $\mathbb{P}^*$
is homogeneous, that the canonical linear extension $<$ is isomorphic to the countable rational order $\eta$. To see this, let $x,y$ be elements of the universal poset such that $x < y$ with respect to the canonical linear extension. Choose elements $x_1, z_1, y_1$ in the universal poset $\mathbb{P}$ such that $x_1 < z_1 <y_1$, and such that the induced  $\mathbb{P}^*$-structure $X_1=((x_1,y_1) \prec, <)$ is isomorphic to 
$X=((x,y) \prec, <)$. Let $\pi_1$ be an isomorphism from $X_1$ to $X$ and use the homogeneity of $\mathbb{P}^*$ to find an automorphism $\pi$ which  extends $\pi_1$. Setting $z=\pi(z_1)$, we have found an element $z$ such that $x < z <y$. In a similar way, we can find, for every element $x$ some elements $y$ and $z$, such that $z <x <y$. It follows that $<$ is a rational total order.

In the sequel, we shall fix a recursive representation of $(\mathbb{P}, <) $.  This means that we can view the underlying set of $\mathbb{P}$ as the natural numbers with both the relations $\prec$ and $ <$ being decidable.  Write $\mathcal{M}(\mathbb{P})$ for the set of linear extensions of $\mathbb{P}$. Clearly $\mathcal{M}(\mathbb{P})$ is a closed subspace of $\mathcal{M}$. 

The automorphism group $\mbox{Aut}(\mathbb{P})$ acts naturally on the space $\mathcal{M}(\mathbb{P})$. Note that the orbit of the canonical linear extension of $\mathbb{P}$ under the action of $\mbox{Aut}(\mathbb{P})$ is contained in the class $\mathcal{Q}$ of rational orders on the natural numbers $\mathbb{N}$. It turns out that all the elements in this orbit are not $\mu$-Martin-L\" of random.
\begin{theorem} 
Fix a recursive representation of the universal poset $\mathbb{P}$ on the natural numbers $\mathbb{N}$. Let $\mathcal{M}(\mathbb{P})$ be the class of linear extensions of $\mathbb{P}$. Write $ML_\mu$ for the set of total orders on $\mathbb{N}$ that are Martin-L\" of random relative to the Glasner-Weiss probability measure $\mu$. Then
\[ML_\mu \cap \mathcal{M}(\mathbb{P}) = \emptyset.\]
\label{th:posetextensions}
\end{theorem}
Proof: Note that for $\xi \in \mathcal{M}$, it is the case that
$$ \xi \in \mathcal{M}(\mathbb{P}) \Longleftrightarrow \forall_{x\in\mathbb{N}}\;\forall_{y\in\mathbb{N}}\; (x \prec y \Rightarrow x <_\xi y),$$
which, since $\prec$ is recursive over $\mathbb{N}$, means that $\mathcal{M}(\mathbb{P})$ is a  $\Pi_1^0$-subset of $\mathcal{M}$. By Proposition \ref{prop:kurtz}, if it were the case that
$ML_\mu \cap \mathcal{M}(\mathbb{P}) \neq \emptyset$,
then $\mu(\mathcal{M}(\mathbb{P}) ) > 0$. However, as has been noted before, Kechris and Soki\v c \cite{KechrisSokic} has recently shown that the automorphism group $\mbox{Aut}(\mathbb{P})$ is {\it not} amenable. Moreover, it is shown by Kechris et al in \cite{Kechrisetal}, that the structural Ramsey theory of finite posets (see, for example, \cite{Fouchposet} ) has the implication that 
this $\mbox{Aut}(\mathbb{P})$-flow on $\mathcal{M}(\mathbb{P})$ is in fact isomorphic to the universal minimal $\mbox{Aut}(\mathbb{P})$-flow; allowing us to infer that $\mu(\mathcal{M}(\mathbb{P}) ) = 0$. This concludes the proof of the theorem.

For a rational order $\eta$, set
\[S_\mu(\eta): =\{\sigma \in S_\infty: \sigma\eta \in ML_\mu\}.\]
By Theorem \ref{th:mlcantor},  if $\eta$ were not rational, the corresponding set $S_\mu(\eta)$ would have been the empty set. 
If $\eta$ is rational, then the class $\mathcal{Q}$ is exactly the orbit of $\eta$ under the action of $S_\infty$. Since both $\mathcal{Q}$ and $ML_\mu$ have $\mu$-measure one, it follows that
$$\mu(\mathcal{Q} \cap ML_\mu)=1,$$
and, therefore, that $S_\mu(\eta) \neq \emptyset$.
 Note that, if $\pi \in S_\infty$, then 
\begin{equation}
 S_\mu(\eta)\pi^{-1} = S_\mu(\pi\eta).
\label{eq:symmetrizers}
\end{equation}
 Indeed, for $\alpha \in S_\mu(\eta)$, we have $\alpha\pi^{-1}(\pi\eta) =\alpha\eta \in ML_\mu$ and hence $\alpha\pi^{-1} \in  S_\mu(\pi\eta)$.  Conversely, if $\tau \in S_\mu(\pi\eta)$, then $\tau\pi\eta \in ML_\mu$, i.e., $\tau \pi \in  S_\mu(\eta) $, and, so, $\tau \in S_\mu(\eta)\pi^{-1}$. 
 
 If $\eta_1,\eta_2 \in \mathcal{Q}$, there is some $\pi \in S_\infty$ such that $\eta_2=\pi\eta_1$. Moreover, if $\eta_1,\eta_2$ 
 were both recursive, the permutation $\pi$ could also be chosen to be recursive. (See Theorem \ref{th:efffraisse}). Write $S_r$ for the class of recursive permutations of $\mathbb{N}$. We let $S_r$ act on the right on the class $\Sigma$ of all sets of the form  $S_\mu(\tau)$ with $\tau$ a recursive rational order on $\mathbb{N}$.  The action is given by
 $$\Sigma \times S_r \longrightarrow \Sigma,$$
 $$(S_\mu(\tau), \pi) \mapsto S_\mu(\tau)\pi^{-1}, \; \pi \in S_r\; \tau \in \mathcal{Q}_r,$$ where  $\mathcal{Q}_r$ denotes the class of all recursive rational orders on $\mathbb{N}$. It follows from the preceding arguments that this $S_r$-action will have a single orbit, i.e, the action is transitive. Set
 $$\mathcal{S} =\bigcup_{\tau \in \mathcal{Q}_r} S_\mu(\tau).$$ If we choose any fixed $\eta \in \mathcal{Q}_r$, we also have
 $$\mathcal{S} =\bigcup_{\pi \in S_r} S_\mu(\eta)\pi^{-1}.$$
 We shall call the permutations in $\mathcal{S}$ {\it Martin-L\" of randomizers}. These are the permutations that transform some recursive rational order to one which is $\mu$-Martin-L\" of random.
 
 If $\pi \in \mbox{Aut}(\mathbb{P})$ and $\eta$ is a recursive rational linear extension of $\mathbb{P}$, then $\pi\eta$ is also a linear extension of $\mathbb{P}$, so, by Theorem \ref{th:posetextensions}, we have $\pi \not\in S_\mu(\eta)$. In particular,
 \begin{theorem}
 No  Martin-L\" of randomizer of a recursive rational linear extension of the universal poset $\mathbb{P}$ can be an automorphism of $\mathbb{P}$.
\end{theorem} 

Let $\mathbb{G}$ be the Fra\" {i}ss\' e  limit of finite graphs. It follows, again  from \cite{Kechrisetal}, that the universal minimal flow of the automorphism group $G$ of $\mathbb{G}$ is isomorphic to the space $\mathcal{M}$ of total orders with group action being the restriction of the $S_\infty$-action to $G$. Of course the Glasner-Weiss measure $\mu$ is also $G$-invariant, which gives another proof of the amenability of $G$. Kechris \cite{Kechris2011} has recently announced that he, in collaboration with Angel and Lyon, have recently shown that $\mu$ is the only $G$-invariant probability measure on $G$. Apparently, the fact that graphs are Hrushovski-structures \cite{Hrushovski} plays an important r\^ole in the proof. It would be interesting to understand which symmetries of the graph $\mathbb{G}$ are in fact Martin-L\" of randomizers.

Write $\mathbb{N}^{[2]}$ for the set of two element subsets of $\mathbb{N}$. We can use $\mathcal{G}=\{0,1\}^{{\mathbb{N}}^{[2]}}$ as a representation of all countable graphs (with the product topology) with underlying set $\mathbb{N}$. The group $S_\infty$ acts continuously on $\mathcal{G}$ , if, for $\pi \in S_\infty$ and $\alpha \in \mathcal{G}$, we set
$$\pi\alpha(\{i,j\})=\alpha (\{\pi^{-1}i,\pi^{-1}j\}),$$ for all $i,j \in \mathbb{N}$. We shall refer to this action as the canonical $S_\infty$-flow on $\mathcal{G}$.

 Write $\mathcal{R} \subset \mathcal{G}$ for the set of $\alpha \in \mathcal{G}$  corresponding to the copies of $\mathbb{G}$, the Fra\" {i}ss\' e limit of finite graphs.  It is well-known that a graph $\alpha \in \mathcal{G} $ defines a copy of $\mathbb{R}$, iff, whenever  $A,B$ are disjoint finite subsets of $\mathbb{N}$, there is some $z$ in the complement of $A\cup B$ such that $\alpha(\{z,a\}) =1$ for all $a \in A$, and $\alpha(\{z,b\})=0$ for all $b \in B$. This result has the following consequence: If $\pi \in S_\infty$ and $\alpha \in \mathcal{R}$, then $\pi\alpha \in \mathcal{R}$. This means that $S_\infty$ acts naturally on the class $\mathcal{R}$ of random graphs. 

Let $\lambda$ be the Lebesque measure on $\mathcal{N}=\{0,1\}^{\mathbb{N}}$ and write $ML_\lambda$ for $\lambda$-Martin-L\" of elements of $\{0,1\}^{\mathbb{N}}$.  Let $\phi$ be any recursive bijection from $\mathbb{N}^{[2]}$  to $\mathbb{N}$.  Then $\phi$ induces a recursive isomorphism $f:=\{0,1\}^\phi$ from  $\mathcal{N}$ to $\mathcal{G}$.  In this way, we find a recursive embedding of $ML_{\lambda}$ into $\mathcal{R}$, the latter being all countable graphs isomorphic to $\mathbb{G}$. (See Fouch\' e and  Potgieter \cite{FouchPot}.) Write $\lambda_1$ for the pushout of $\lambda$ under $f$. Then $f$ induces a recursive isomorphism between $ML_\lambda$ and $ML_{\lambda_1}$. These are all the random graphs encoded via $\phi$ by some infinite binary string which is $\lambda$-Martin-L\"of random. Note that $\lambda_1$ is a $S_\infty$-invariant measure under the canonical $S_\infty$-flow .

In recent work, Petrov and Vershik (2010) \cite{PetrovVershik1} have made a systematic study of the invariant measures relative  the canonical $S_\infty$-flow restricted to the class $\mathcal{R}$. They have identified a continuum of $S_\infty$-invariant measures on $\mathcal{R}$ . The study of Martin-L\" of randomness and the associated randomizers relative to the computable $S_\infty$-invariant measures,  is at present under investigation by the author.

\subsection*{Acknowledgements} 
The research is based upon work supported by the National Research Foundation (NRF) of South Africa. Any opinion, findings and conclusions or recommendations expressed in this material are those of the author and therefore the NRF does not accept any liability in regard thereto.


\end{document}